\documentclass[final,a4paper,5p]{elsarticle}
\pdfoutput=1
\usepackage[colorlinks=true,linkcolor=blue,citecolor=blue,urlcolor=blue]{hyperref}
\usepackage[utf8]{inputenc}
\usepackage[english]{babel}
\usepackage{amsmath}
\usepackage{amsfonts}
\usepackage{amssymb}
\usepackage{graphicx}
\usepackage{epstopdf}
\usepackage{listings}
\usepackage{empheq}
\usepackage{url}

\DeclareUnicodeCharacter{00B5}{\ensuremath{\mu}}
\usepackage{color}
\usepackage{ifthen}
\newsavebox{\picturebox}
\newlength{\pictureboxwidth}
\newlength{\pictureboxheight}
\newcommand{\includeimage}[1]{
	\savebox{\picturebox}{\includegraphics{#1}}
	\settoheight{\pictureboxheight}{\usebox{\picturebox}}
	\settowidth{\pictureboxwidth}{\usebox{\picturebox}}
	\ifthenelse{\lengthtest{\pictureboxwidth > .95\linewidth}}
	{
		\includegraphics[width=.95\linewidth,height=.80\textheight,keepaspectratio]{#1}
	}
	{
		\ifthenelse{\lengthtest{\pictureboxheight>.80\textheight}}
		{
			\includegraphics[width=.95\linewidth,height=.80\textheight,keepaspectratio]{#1}
			
		}
		{
			\includegraphics{#1}
		}
	}
}

\usepackage{animate} % This package is required because the wxMaxima configuration option
\definecolor{labelcolor}{RGB}{100,0,0}
\usepackage{fancyhdr}
\begin{document}
\renewcommand{\today}{December 21\textsuperscript{st}, 2018}

\begin{abstract}
A novel method rooted in the classical Schwarz-Christoffel transformation from the disk is introduced, which allows for fast and accurate solution of potential field problems in possibly inhomogeneous and multiply connected domains: this is for sure its most outstanding feature, circumventing the barriers that have increasingly restricted the scope of conformal mappings in applications since the advent of computers and purely numerical methods. An example problem, derived from a case of practical interest, is analyzed and results are compared with those obtained from FEA.
\end{abstract}


\begin{frontmatter}
	\title{An efficient procedure for solving potential field problems: \\ the Conformal Boundary Differences Method}
	
	\author{Stefano Costa\fnref{fn1}}
	\ead{stefano.costa@ieee.org}
	
	\fntext[fn1]{IEEE, member; Piacenza 29122 Italy. {\includegraphics[scale=0.5,keepaspectratio]{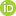}} \url{ https://orcid.org/0000-0002-9434-356X}}
	
	\journal{arXiv}
\end{frontmatter}

\section{Introduction}
Conformal mapping has a long and successful history as a mathematical tool for solving field problems in physics and engineering. It is by far superior to any other computational method whenever it can provide an \emph{analytical} solution, showing explicitly how it changes along with its variables. A closer look to the reasons that have gradually restricted its usage in applications will serve as a foreword to the motivation for the present work.

The trouble of determining the accessory quantities (the "parameter problem") has been nowadays superseded by modern computers and the development of dedicated software packages, the most consistent and popular being undoubtedly the Schwarz-Christoffel Toolbox for \textsc{Matlab} by T. A. Driscoll \cite{SCTB}, which has made these kind of mappings a matter of a few mouse clicks. Still some difficulties can arise with the \emph{crowding of the preverties} phenomenon: this is met either reverting to a more suitable canonical domain, or by means of the CRDT algorithm by Driscoll and Vavasis \cite{DV98}; regrettably, the latter is not available for unbounded regions.

Multiply connected domains have always represented an issue, unless some sort of symmetry could be exploited, that renders the reduced domain simply connected. This was, until recent theoretical developments \cite{DEP04} and the celebrated breakthroughs of D. Crowdy in particular \cite{C05}, \cite{C07}, which provide explicit formulas for such situations. Whereas these represent a huge deepening in theoretical understanding, their translation to effective computer algorithms are all but trivial and, to our best knowledge, nonexistent at least as packages like the aforementioned SC-Toolbox, although one was developed in FORTRAN 77 for doubly connected regions \cite{HU98}.

Finally, the great wall of nonuniform domains: it is nearly impossible to deal with them analytically, except in very simple cases, and one must compromise a big deal in order to force certain assumptions. Among the three, this is the biggest drawback that makes nowadays conformal mappings a niche method in applications: the textbook by Schinzinger and Laura \cite{SL03}, a collection of advanced techniques for a broad range of applications, illustrates the treatment of nonuniform media by a few particular examples only, often requiring beforehand approximations.

On the engineering side, we refer to the work of a research group in Pavia (Italy) that, since the early 2000s, has put great effort and creativity in the analysis of electric and magnetic fields by means of numerical Schwarz-Christoffel transformations. Beginning with simple yet innovative applications of the SC-Toolbox \cite{CCBS00}, they came to deal with doubly connected \cite{CBS09} and inhomgeneous domains \cite{CBMS10}, \cite{CB17}. As remarkable as this progress is, here too we meet the distinctive limitations discussed above: doubly connected SC transformations map to the annulus only, thus requiring the analyst a notable degree of skill in order to cascade other mappings and manipulations; furthermore, modeling different media by means of current sheets is not always feasible, or desirable. In the Authors' words "\emph{computational procedures have to be developed \emph{ad hoc} for a given class of problems}", and "\emph{a severe drawback is that, in principle, the transformation technique requires to model homogeneous materials}" \cite{CBMS10} p. 66. A common situation is the use of Schwarz-Christoffel transformations limited to ancillary operations prior to analyses via FDM or FEM in domains inhomogeneous or having various types of boundary conditions, where the last can work smoothly.
\\

The aim of the present work is a standardized procedure for the solution of potential field problems by means of Schwarz-Christoffel transformations, capable of circumventing the severe restrictions discussed above. We adopt a constructive approach, starting from a  model problem, to focus tightly on its steps and put it to the test; therefore we drop the treatment of field sources and forcibly restrict ourselves to Laplacian fields:
\begin{subequations}
	\begin{align}
 	& \text{Laplace's equation:} \label{eqn:lapeq} \ \nabla^2 u(z)=0,\, z\in\mathbb{P} \\
	& \text{Dirichlet boundary condition:} \ u(z) = \psi_D(z),\, z\in\Gamma_D \\
	& \text{Neumann b.c.:} \label{eqn:neucond} \ \dfrac{\partial u}{\partial n}\equiv\nabla u(z)\cdot\hat{n}=\gamma(z),\, z\in\Gamma_N
	\end{align}
\end{subequations}
where $u(z)$ is the potential to be determined; $\mathbb{P}$ is the domain problem, possibly inhomogeneous and multiply connected; $\hat{n}$ is the outward normal to the boundary $\partial\mathbb{P}$; $\Gamma_D$ and $\Gamma_N \subseteq \partial\mathbb{P}$.

The rationale is: we want these transformations to work in situations they are considered unsuited to, so we must be willing to some trade-offs, rethinking the way they have always been employed. Accordingly, we decide to weaken the representation of a domain no longer seen as a whole, characteristic of conformal mappings, and replace overall continuity with some targeted discretization, for a broader and less case dependent applicability.

\section{The Model Problem}
\begin{figure}[!t]
	\centering
	\includegraphics[width=1\linewidth]{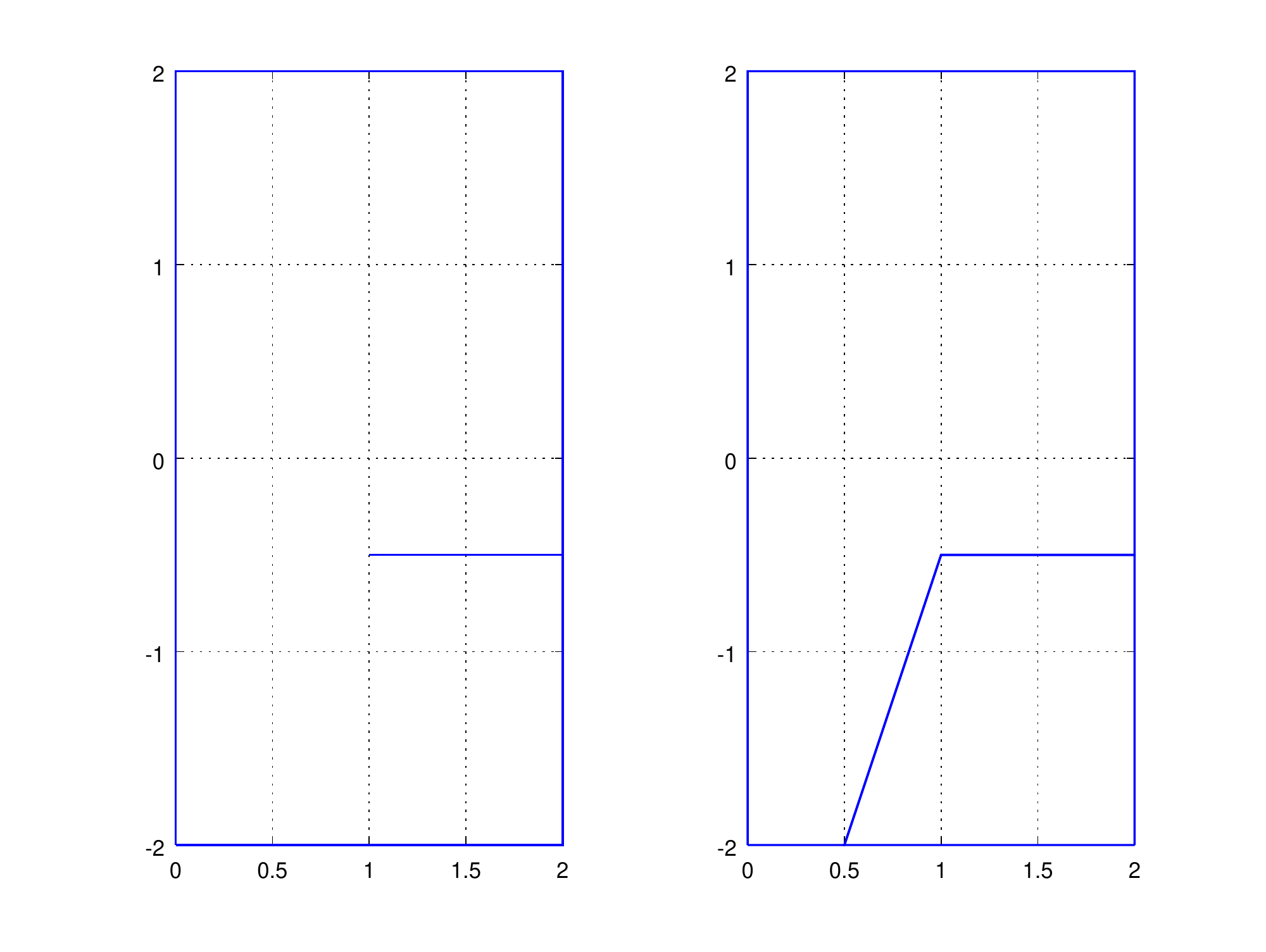}
	\caption{The thin microstrip (horizontal), surrounded by dielectric and shielded by an equipotential conductor (lengths in arbitrary units); the vertical right side of the rectangle is a line of symmetry with homogeneous Neumann boundary condition $\partial u/\partial n=0$. In the picture on the left the dielectric is necessarily homogeneous, whereas on the right it may be not, the diagonal line representing the interface between possibly different media.}
	\label{fig:scgeom}
\end{figure}

Our testbed is taken from \cite{CB17}, with some variants\footnote{In the referenced article, the microstrip lies on a dielectric support of rectangular shape. Our variant, whereas not really representative of real devices, shows how diagonal (arbitrarily directed) boundaries are easily dealt with.} helpful to better illustrate the characteristics of our solution method, as depicted in figure \ref{fig:scgeom}.
We want to calculate the capacitance of a thin microstrip (zero thickness) surrounded by dielectric, possibly inhomogeneous, and shielded by an equipotential conductor. Only half of the object needs to be analyzed by exploiting its vertical symmetry, as appears from figure \ref{fig:field_all}, therefore we impose an homogeneous Neumann boundary condition along the whole right boundary of the rectangle. Additionally, we want to find how the electric potential varies along the vertical line with coordinate $x=0.999$ arbitrary units: this is extremely close to the end of the microstrip with $x=1$, where the field gradient becomes singular, and needs particular care to be handled correctly. The rectangle is the problem domain $\mathbb{P}$, possibly split into two subdomains $\mathbb{A}$ (upper, bigger) and $\mathbb{B}$ (lower, smaller); the strip and the shield are boundaries of the Dirichlet type $\Gamma_D$, whereas the line of symmetry and the interface between subdomains are of the Neumann type $\Gamma_N$.

\begin{figure}[!t]
	\centering
	\includegraphics[width=0.9\linewidth]{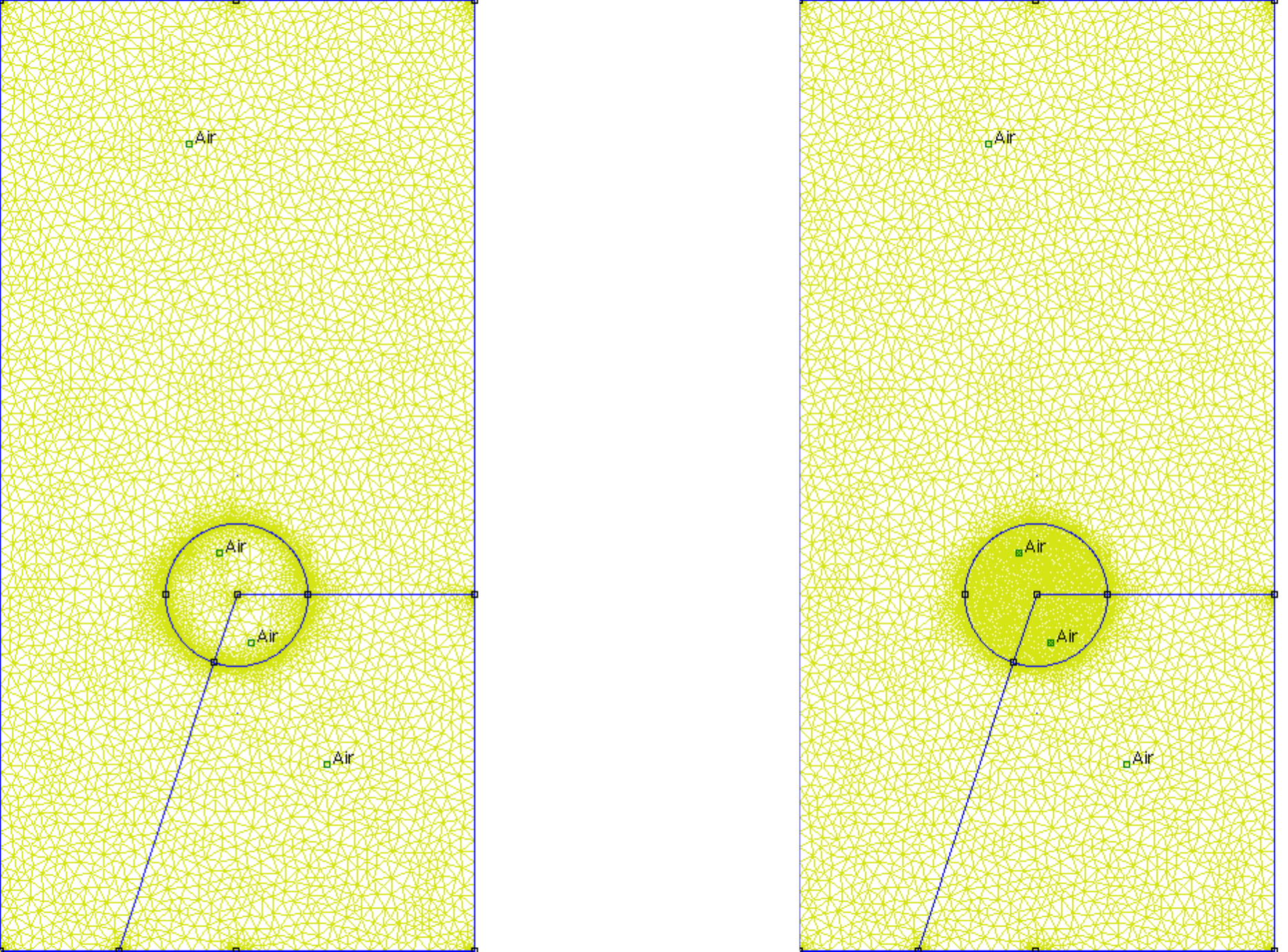}
	\caption{On the left, the problem domain meshed by FEMM with 6357 nodes/12436 elements (automatic generation); on the right, with 17306 nodes/34334 elements (element size of 0.005 arbitrary units forced at the strip end).}
	\label{fig:FEMM_GRID}
\end{figure}

Figure \ref{fig:FEMM_GRID} shows the problem domain as meshed by the FEA software FEMM \cite{FEMM}, with automatic generation and when imposing an element size of 0.005 units in the circular region of radius 0.2 units centered at the strip end. Forcing the element size causes a huge increase in the number of nodes, thus allowing better understanding of FEM behaviour near critical points.

\section{Foundations of the Method}
In what follows we move at expedite pace through the elements of potential theory, conformal mappings and discretization of PDEs; we refer to the work by Binns and Lawrenson \cite{BL73} for a coverage of all these topics within a single book. Also, as the field of complex numbers $\mathbb{C}$ is isomorphic to the vector space $\mathbb{R}^2$, we freely switch between the notations $z=x+jy$, $z=\rho e^{j\theta}$ and $z=(z_x,z_y)$ as needed by the discussion. 

To begin with, if our problem is well posed we know from potential theory that a solution $u(z)=\psi(z)$ exists and is unique also on $\partial\mathbb{P}$; but this also means that (\ref{eqn:neucond}) can be restated yet as another Dirichlet b.c. once we have determined the distribution $\psi_N(z)$:
$$u(z) = \psi_N(z),\, z\in\Gamma_N \quad \text{also satisfying }\dfrac{\partial\psi_N}{\partial n}=\gamma(z)$$
This simple fact steers our attention to the unit disk $\mathcal{D}:|t|=|u+jv|<1$, where it is possible to obtain a solution to any problem of the Dirichlet type by means of the \emph{complex potential function of Schwarz}\footnote{One might also recall the \emph{Poisson integral} here, restricted to the real potential $\psi$: $$\psi(t)=\psi(\rho e^{j\beta})=\dfrac{1}{2\pi}\int_0^{2\pi}\frac{1-\rho^2}{1-2\rho\cos(\beta-\theta)+\rho^2}\,\psi(\theta)\,d\theta$$ but the reason for preferring the former will be apparent later.}. The formula is truly notable, being capable of determining the potential at \emph{any} interior point $t$ when it is known along the boundary alone:
\begin{align}
\label{eqn:schwarzian} \psi(t)+j\phi(t) = \dfrac{1}{2\pi} \int_0^{2\pi}\frac{e^{j\theta}+t}{e^{j\theta}-t}\,\psi(\theta)\,d\theta
\end{align}
where $\phi(t)$ is the \emph{flux function}, harmonic conjugate of the potential, and $\psi(\theta)$ is the potential distribution along $\partial\mathcal{D}$.

\begin{figure}[!t]
	\centering
	\includegraphics[width=1\linewidth]{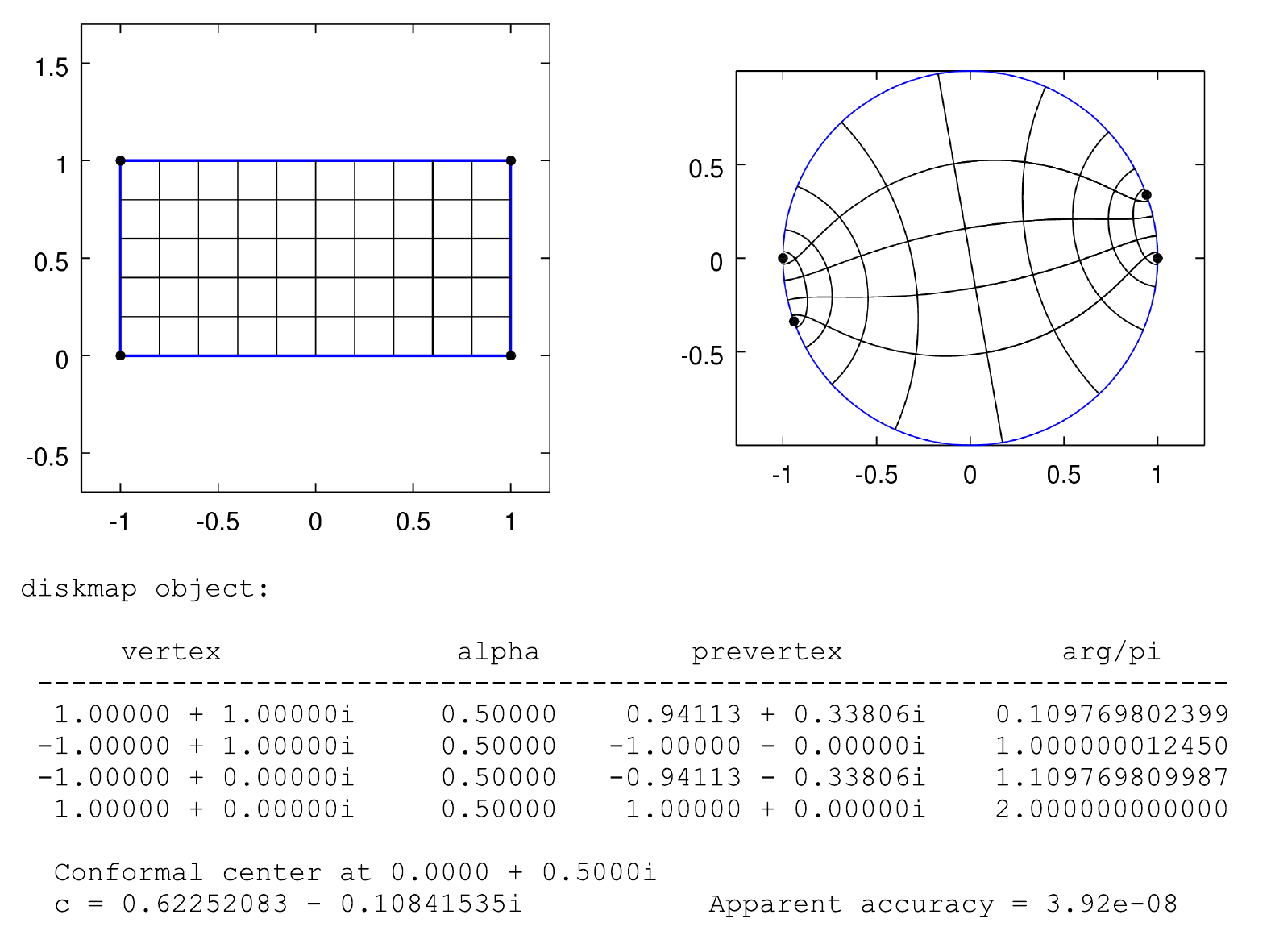}
	\caption{Example of conformal mapping from the unit disk $\mathcal{D}$ to the rectangle $\mathbb{P}$ obtained by the SC-Toolbox. \texttt{alpha} represents the $\{\alpha_k\}$ in (\ref{eqn:scdisk}), i.e. the turning angles at the vertices of the polygon. Both vertices and prevertices are marked with dots: each prevertex $t_k\in\partial\mathcal{D}$ maps to a vertex $z_k=F(t_k) \in \partial\mathbb{P}$, and the same do the curvilinear streams, representing flux and equipotential lines, to the orthogonal grid. Here one can see clearly that orthogonality between them is maintained after transplantation from (or to) the disk, and this makes Schwarz-Christoffel mapping a formidable tool for dealing with Laplace's equation.}
	\label{fig:scgrid}
\end{figure}

The Schwarz-Christoffel mapping from the unit disk is the next tool that fits into the box naturally:
\begin{align}
\nonumber z & = F(t) = z_0+C\int_{t_0}^t \prod_{k=1}^n \left(t_k-\tau\right)^{\alpha_k-1}\,d\tau = \\
& = z_0+\int_{t_0}^t f(\tau)\,d\tau \label{eqn:scdisk}
\end{align}
Suffice to recall here that the formula above maps \emph{conformally} the unit disk $\mathcal{D}$ in the $t$-plane to any simply connected polygon $\mathbb{P}$ in the $z$-plane. Figure \ref{fig:scgrid} may serve as a (tiny) refresher for the reader familiar with this subject; and we refer to \cite{DT02} anyone willing to delve deeper (and broader).

\begin{figure}[!t]
	\centering
	\includegraphics[width=0.9\linewidth]{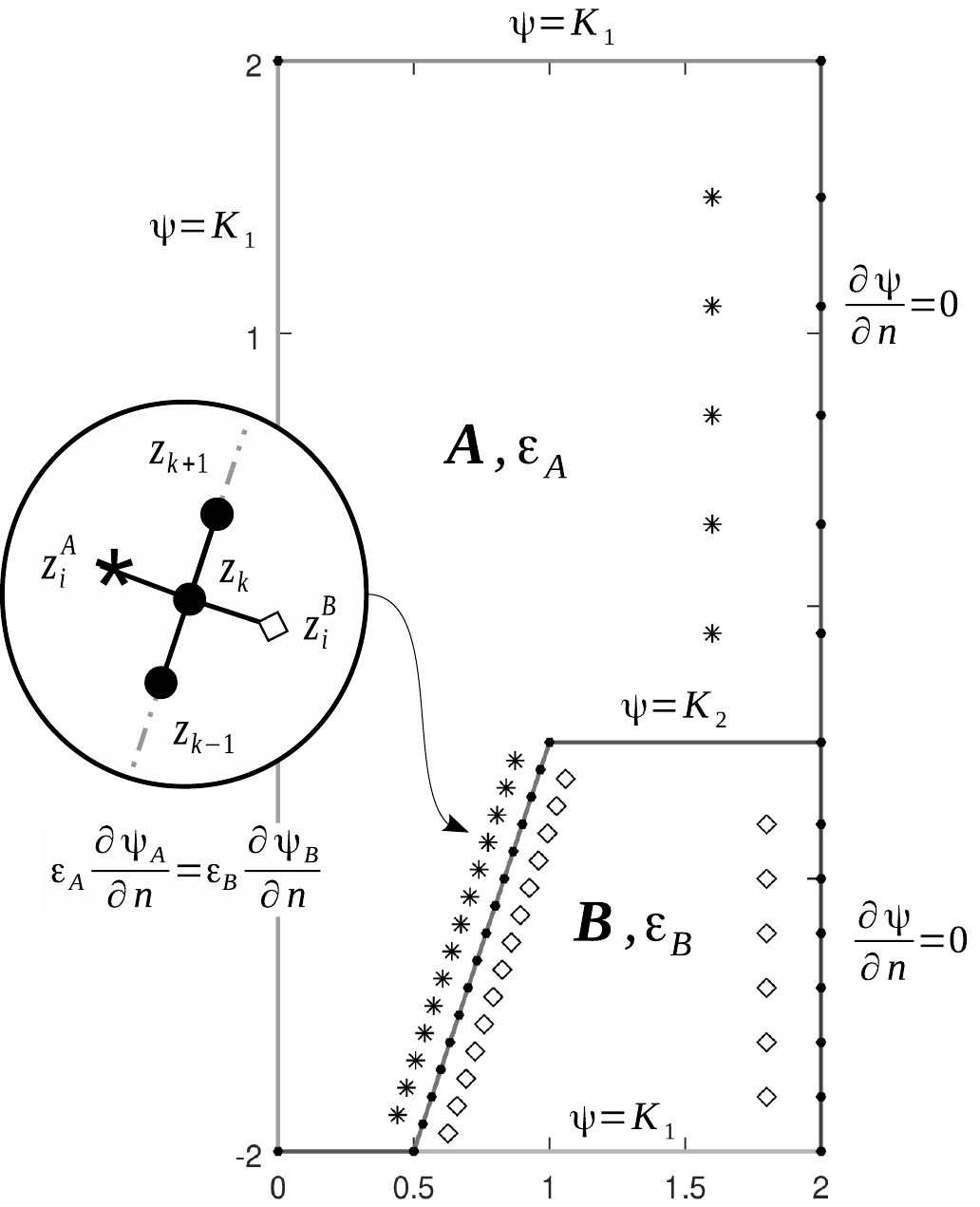}
	\caption{The domain $\mathbb{P}$ is partitioned into two simply connected polygons $\mathbb{A}$ and $\mathbb{B}$. The upper, left and lower sides, and the horizontal strip, are at a defined potential, i.e. are boundaries of the Dirichlet type $\Gamma_D$, whereas the sides on the right, and the interface, are of the Neumann type $\Gamma_N$. Discretization and FD are used on $\Gamma_N$ only to enforce conditions on $\partial\psi/\partial n$: the $\{z_k\}$ are represented by dots, and the $\{z_i\}$ by asterisks and diamonds in $\mathbb{A}$ and $\mathbb{B}$ respectively. Stencils of different sizes and possibly asymmetric can be employed with no difficulties arising from mesh generation; note also that dealing with diagonal boundaries comes at no price in terms of additional machinery.}
	\label{fig:schematic}
\end{figure}

$F(t)$ is a bijective complex valued function of a complex variable, \emph{holomorphic} and \emph{conformal}, i.e. analytic and capable of preserving angles between intersecting lines after transformation. These facts when applied to the \emph{complex potential field}\footnote{The variable is left unspecified intentionally.} $w=\psi+j\phi$ lead to the following important results amongst many others (see e.g. \cite{H74}):
\begin{subequations}
	\begin{align}
	& \text{Invariance of potential:}  \label{eqn:pot} \ w(z) = w(F(t)) = w(t) \\
	& \text{Invariance of Laplace's eqn.:} \label{eqn:lap} \ \nabla^2\psi(z) = \nabla^2\psi(t) = 0 \\
	& \text{Law of potential gradients:} \label{eqn:grad} \ \nabla\psi(z) \cdot \overline{f(t)}= \nabla\psi(t)
	\end{align}
\end{subequations}
We now describe how to make effective use of formulas (\ref{eqn:schwarzian})-(\ref{eqn:grad}) before seeing them in action in the next section: our goal is a procedure to determine the unknown $\psi_N(z)$. Let's consider the situation sketched in figure \ref{fig:schematic}: the idea is that we can enforce constraints at all boundaries of the Neumann type $\Gamma_N$, namely borders and interfaces between subdomains, by a suitable discretization and applying Finite Differences there.
One can start with an arbitrary distribution of potential at $\{z_k\}\in\Gamma_N$ and  $\{z_i\}\in\mathbb{P}$, say all 0's, and converge to a solution via successive over-relaxation (SOR) iterations. The good news is, we don't really need a grid here to find and update the values of $\psi$. At each step, we use the images $t_i=F^{-1}(z_i)$ in the unit disk $\mathcal{D}$ and update the values $\psi(t_i)=\psi(z_i)$ by using (\ref{eqn:schwarzian}) with a distribution $\psi(\theta)$ somehow rebuilt from the current $\psi(t_k)=\psi(z_k)$, where $t_k=F^{-1}(z_k)\in\partial\mathcal{D}$. Then go back to $\mathbb{P}$, update $\psi(z_k)$ via FD, and repeat. This is the overall picture; let's go through the steps of the algorithm in detail:
\begin{enumerate}%[\IEEEsetlabelwidth{9)}]
	\item Define a set of disjoint polygons  $\{\mathbb{P}_m\}:\bigcup_m\mathbb{P}_m=\mathbb{P}$, each having Dirichlet or Neumann boundaries, or an arbitrary mixture of the types.
	\item Discretize each boundary of the Neumann type, including interfaces between polygons, with suitable steps, the smaller the better, by positioning sets of points $\{z_k\}_m$, and corresponding $\{z_i\}_m$ orthogonally with respect to the direction of the boundary. Note that orthogonality is required \emph{locally}\footnote{Clearly (\ref{eqn:lap}) must hold irrespective of the particular orientation.}, as opposite to classical FD schemes where diagonal or irregular boundaries not matching the mesh grid often represent a source of coarse approximation; also, different stencil dimensions, possibly asymmetric, can be used where needed without the inconvenience of having to "bridge" them.
	\item For each $\mathbb{P}_m$ compute a transformation $F_m(t)$ from the unit disk $\mathcal{D}_m$; polygons must be such shaped, that numerical Schwarz-Christoffel mappings are not affected by the crowding phenomenon. Thus we can exploit (\ref{eqn:pot}) and transfer our problem to the disk.
	\item\label{point:stepA} Do not apply (\ref{eqn:schwarzian}) as is: it would bring unnecessary numerical difficulties only. Instead by noting that for $t=0$ it holds
	\begin{align*}
	%\label{eqn:psi0}
	\psi(0)=\dfrac{1}{2\pi}\int_0^{2\pi}\psi(\theta)\,d\theta
	\end{align*}
	that is a much simpler integral of a real valued function, for each $t_i\in\mathcal{D}_m$ do remap the $\{t_k\}_m\in\partial\mathcal{D}_m$ via a M\"{o}bius transformation\footnote{Yet another conformal mapping, featuring the same relevant properties recalled before.} that leaves their sequence in $\partial\mathcal{D}_m$ unaltered while sending $t_i$ to the origin:
	\begin{align*}
	& t'=e^{j\theta'}=g_i(t)=\dfrac{t-t_i}{1-\overline{t_i}\,t} \\
	& \{t_k\}_m\mapsto\{t'_k\}_m=\{e^{j\theta'_k}\}_m
	\end{align*}
	hence:
	\begin{align}
	\nonumber \psi(t_i) & =\psi(g_i(t_i))=\psi(0)= \\
	& =\dfrac{1}{2\pi}\int_0^{2\pi}\psi(\theta')\,d\theta'=\psi(z_i) \label{eqn:itert}
	\end{align}
	Otherwise stated, in this step build collections of M\"{o}bius mappings $\{g_i(t)\}_m$ for a much more accurate evaluation of (\ref{eqn:schwarzian}).
	\item\label{point:stepB} Update the values of potential at $\Gamma_N$ through the well-known FD algorithm, here given for 5-points regular stencils:
	\begin{align}
	&\nonumber \psi(z_{k-1})+\psi(z_{k+1})+\ldots & \\
	&\ldots+c_A\psi(z_i^A)+c_B\psi(z_i^B)-4\psi(z_k)=0 & \label{eqn:FD}
	\end{align}
	where $A$ and $B$ indicate the two polygons at the interface, and coefficients $c_A$ and $c_B$ account for a possible difference in the media\footnote{In the case of two dielectrics with permittivity $\varepsilon_A$ and $\varepsilon_B$, $c_A=2\varepsilon_A/(\varepsilon_A\varepsilon_B)$ and $c_B=2\varepsilon_B/(\varepsilon_A\varepsilon_B)$.}; clearly if the regions are uniform $c_A=c_B=1$, and for a homogeneous Neumann b.c. in particular it holds:
	$$ \psi(z_{k-1})+\psi(z_{k+1})+2\psi(z_i)-4\psi(z_k)=0 $$
	\item Repeat step \ref{point:stepA}, equation (\ref{eqn:itert}) and step \ref{point:stepB}, equation (\ref{eqn:FD}) until the desired accuracy on $\psi_N(z)$ is reached; more on this point later.
\end{enumerate}
We have come up with an algorithm that exploits the purely geometric properties of conformal mappings, in order to provide a solution via FD on (sub)domain boundaries only: hence we give it the name of Conformal Boundary Differences Method, the \textbf{CBDM} from now on.
\\

Last but not least, some words on equation (\ref{eqn:grad}). For a given transformation from the unit disk, once the prevertices have been determined, the computation of $f(t)=dz/dt$ is immediate, and so is the gradient in the $z$-plane once we have determined it in the $t$-plane.

Suppose that, by some more or less refined method, we have determined a distribution $\psi(\theta)$. The mean value of potential $\langle\psi_n\rangle$ on any interval $\Delta\theta_n=\theta_{n+1}-\theta_n$ satisfies:
\begin{align*}
\int_{\theta_n}^{\theta_{n+1}}\psi(\theta)\,d\theta = \langle\psi_n\rangle\Delta\theta_n
\end{align*}
Integrating (\ref{eqn:schwarzian}) for equipotential sections gives:
\begin{align*}
w(t) = \dfrac{-1}{2\pi}\sum_n\Delta\theta_n\langle\psi_n\rangle - \dfrac{j}{\pi}\sum_n\langle\psi_n\rangle \log\dfrac{e^{j\theta}-t}{e^{j\theta}-t}
\end{align*}
The gradient of $\psi=\Re(w)$ in the origin cannot be determined in polar coordinates: so we replace $t$ with $u+jv$, feed it to Maxima \cite{MAXIMA} (see Appendix) and finally get:
\begin{align}
& \nonumber \nabla\psi(0) =\nabla\psi(u+jv)|_0 = \left(\partial_u \psi(t)\,, \partial_v \psi(t)\right)|_0 = \\
& =\sum_n \frac{\langle\psi_n\rangle}{\pi} \left( \sin(\theta_{n+1})-\sin(\theta_n)\,,\cos(\theta_n)-\cos(\theta_{n+1}) \right) \label{eqn:gradcbdm}
\end{align}
It seems that such an easy formula for field gradient can not be worked out if one considers the Schwarz-Christoffel mapping from the upper half plane in place of that from the unit disk; also, it isn't clear at all how to deal with $\psi(\theta)$ correctly, even less with $\int\psi\,d\theta$, in presence of a prevertex at infinity: and this is our ultimate reason for choosing the disk. Whereas derived from piecewise average of potential, (\ref{eqn:gradcbdm}) avoids numerical differentiation and we have found it precise to an excellent degree.

\section{Application of the Algorithm}
Before passing on to numerical results, there are at least two topics we must touch upon. The first is about the evaluation of $\int\psi\,d\theta$ in step \ref{point:stepA} of the algorithm. Between any two points $z_0$ and $z_1$  of a stencil it is natural to consider a linear variation of potential:
\begin{equation}
\label{eqn:linpot}\psi(z) = \dfrac{\psi(z_1)-\psi(z_0)}{z_1-z_0}(z-z_0)+\psi(z_0)
\end{equation}
Let now be the following: $\psi(z_0)=\psi_0$, $\psi(z_1)=\psi_1$, $\psi_1-\psi_0=\Delta\psi$ and $z_1-z_0=\Delta z$; from (\ref{eqn:scdisk}) it follows:
\begin{align*}
(z-z_0) & =\int_{t_0}^t f(\tau)\,d\tau \\
\psi(z) & =\dfrac{\Delta\psi}{\Delta z}\int_{t_0}^t f(\tau)\,d\tau + \psi_0 = \psi(t)
\end{align*}
The last integral tells us that $\psi(t)$ is \emph{not} linear with $t$, therefore an analogous to (\ref{eqn:linpot}) on the disk may be used bearing it in mind that it's a coarse (yet effective, as we will see) approximation. As we keep moving along the disk boundary it is $t=e^{j\theta}$ and:
\begin{align*}
\Delta\theta & =(\theta_1-\theta_0) \\
(z-z_0) & =\int_{\theta_0}^{\theta} f(e^{j\theta'})\,d\theta' \approx \dfrac{\Delta z}{\Delta\theta}(\theta-\theta_0) \\
\psi(\theta) &  \approx\dfrac{\Delta\psi}{\Delta z}\left(\dfrac{\Delta z}{\Delta\theta}(\theta-\theta_0)\right) +\psi_0 = \dfrac{\Delta\psi}{\Delta\theta}(\theta-\theta_0)+\psi_0
\end{align*}
hence, omitting some tedious manipulations\footnote{The M\"{o}bius transformation in between, acting as a mere rearrangement of points, doesn't quite change the final result in its form but in the width of $\Delta\theta$, so we omit it in order to keep notations light.}:
\begin{align}
\nonumber \int_{\theta_0}^{\theta_1}\psi(\theta)\,d\theta & =\dfrac{\Delta\psi}{\Delta\theta}\int_{\theta_0}^{\theta_1}(\theta-\theta_0)\,d\theta + \int_{\theta_0}^{\theta_1}\psi_0\,d\theta= \\
& =\dfrac{\psi_1+\psi_0}{2}\,\Delta\theta \label{eqn:potapprox}
\end{align}
As expected, the integral of the linear-varying approximated potential between two points equals their mean value times their angular distance along the circumference. The potential at the center of the disk is finally given by summing (\ref{eqn:potapprox}) for all intervals $[\theta_n,\theta_{n+1}]$ partitioning $\partial\mathcal{D}$:
\begin{align}
\nonumber \sum_n\Delta\theta_n & =\sum_n(\theta_{n+1}-\theta_n) = 2\pi \\
\nonumber \psi(0) & =\dfrac{1}{2\pi}\sum_n\int_{\theta_n}^{\theta_{n+1}}\psi(\theta)\,d\theta= \\
& \label{eqn:potcbdm}  =\dfrac{1}{2\pi}\sum_n\left(\dfrac{\psi_{n+1}+\psi_n}{2}\,\Delta\theta_n\right)
\end{align}

The second point is about our implementation of SOR. Our main source of ideas for numerical \emph{experiments} have been the Numerical Recipes \cite{NR07}, and we have found that the following leads to a large save in number of iterations:
\begin{align*}
& \text{Size of the problem:} \ J= \text{nr. of points on } \Gamma_N \\
& \text{Jacobi radius:} \ \rho_J = 0.999\left(1-\dfrac{\pi^2}{2J^2}\right) \\
& \text{SOR parameter, (optimal) limit value:} \ \omega_{lim} = \dfrac{2}{1+\pi/J}
\end{align*}
With the proposed $\rho_J$ the resulting value of $\omega_{lim}$ is pretty "aggressive", i.e. close to the upper limit value of 2, yet we have never met any problem with stability. We update potential value at the nodes $\{z_k\}$ alternating odd and even values of $k$, with Chebychev acceleration of the SOR parameter $\omega$:
\begin{align*}
& \omega_{odd}^{(0)}=1 \\
& \omega_{even}^{(0)}=1/(1-0.5\,\rho_J^2) \\
& \omega_{odd}^{(n+1)}=1/(1-0.25\,\omega^{(n)}_{even}\,\rho_J^2) \\
& \omega_{even}^{(n+1)}=1/(1-0.25\,\omega_{odd}^{(n)}\,\rho_J^2) \\
& \vdots \\
& \omega^{(\infty)} = \omega_{lim}
\end{align*}
At each step, we calculate the residual $\xi^{(n)}(z_k)$ at each node and their summation $\sum_k\xi$; we stop iterating when maximum residual and $\sum_k\xi$ get as small as required at the same time. The update of potential for a symmetrical stencil is as usual:
\begin{align*}
\psi^{(n+1)}(z_k) = \psi^{(n)}(z_k)+\omega^{(n)}\dfrac{\xi^{(n)}(z_k)}{4}
\end{align*}
for odd and even values of $k$.
\\

In all the subsequent analyses, the shield is at potential $\psi=0$ V and the strip at $\psi=1$ V. We require a tolerance on numerical disk mappings within $10^{-9}$, and on CBDM SOR residuals within $10^{-6}$. Our code is written for Octave \cite{OCTAVE}, which can run the SC-Toolbox ver. 2.1 after some tweaks of little effort.
\begin{figure*}[!t]
	\centering
	\includegraphics[width=0.8\linewidth]{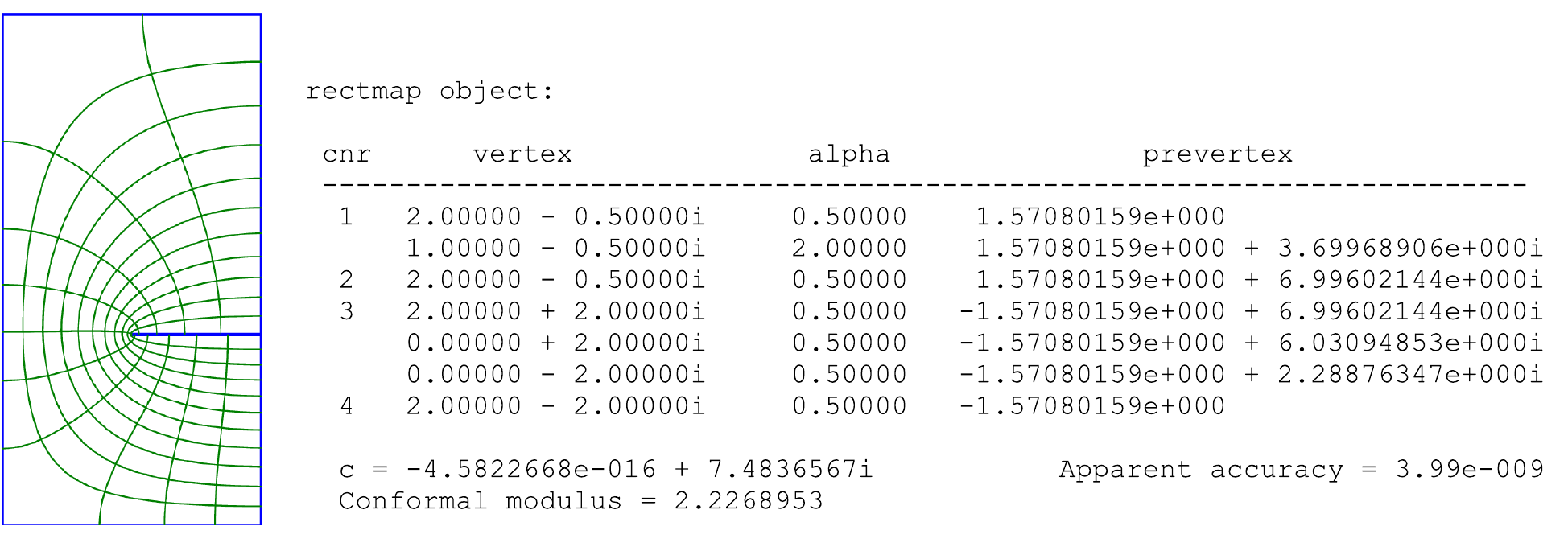}
	\caption{The map from the homogeneous rectangle domain obtained by the SC-Toolbox; equipotential and flux lines are shown.}
	\label{fig:screct}
\end{figure*}
Our first comparison is carried out for a domain $\mathbb{P}$ homogeneous  ($\varepsilon_A=\varepsilon_B=\varepsilon_0$) and, as regards the CBDM, considered as a whole as in figure \ref{fig:scgeom}-left. The values of potential along the vertical line with coordinate $x=0.999$ units, obtained from the methods in table \ref{tab:1}, are compared\footnote{Here we don't alternate odd-even nodes nor use Chebychev acceleration of the SOR parameter.} with those obtained from a conformal transformation from the rectangle in figure \ref{fig:screct}, requiring no additional manipulations and for this reason being considered "exact".
\begin{table}[h]
	\centering
	\begin{tabular}{|c|c|p{4cm}|} \hline
		\textbf{Label}	& \textbf{Method}  & \textbf{Details}  \\ \hline\hline
		FEMM3k	& FEMM & coarser mesh, see fig. \ref{fig:FEMM_GRID}-left \\ 
		FEMM17k	& FEMM & finer mesh, see fig. \ref{fig:FEMM_GRID}-right \\ 
		CBDM005	& CBDM & step of 0.05 units on $\Gamma_N$; 129 SOR iterations \\
		CBDM002	& CBDM & step of 0.02 units on $\Gamma_N$; 333 SOR iterations \\ \hline
	\end{tabular}
	\caption{Description of Methods in Figure \ref{fig:hom1-V}}
	\label{tab:1}
\end{table}
\begin{figure*}[!t]
	\centering
	\includegraphics[width=0.8\linewidth]{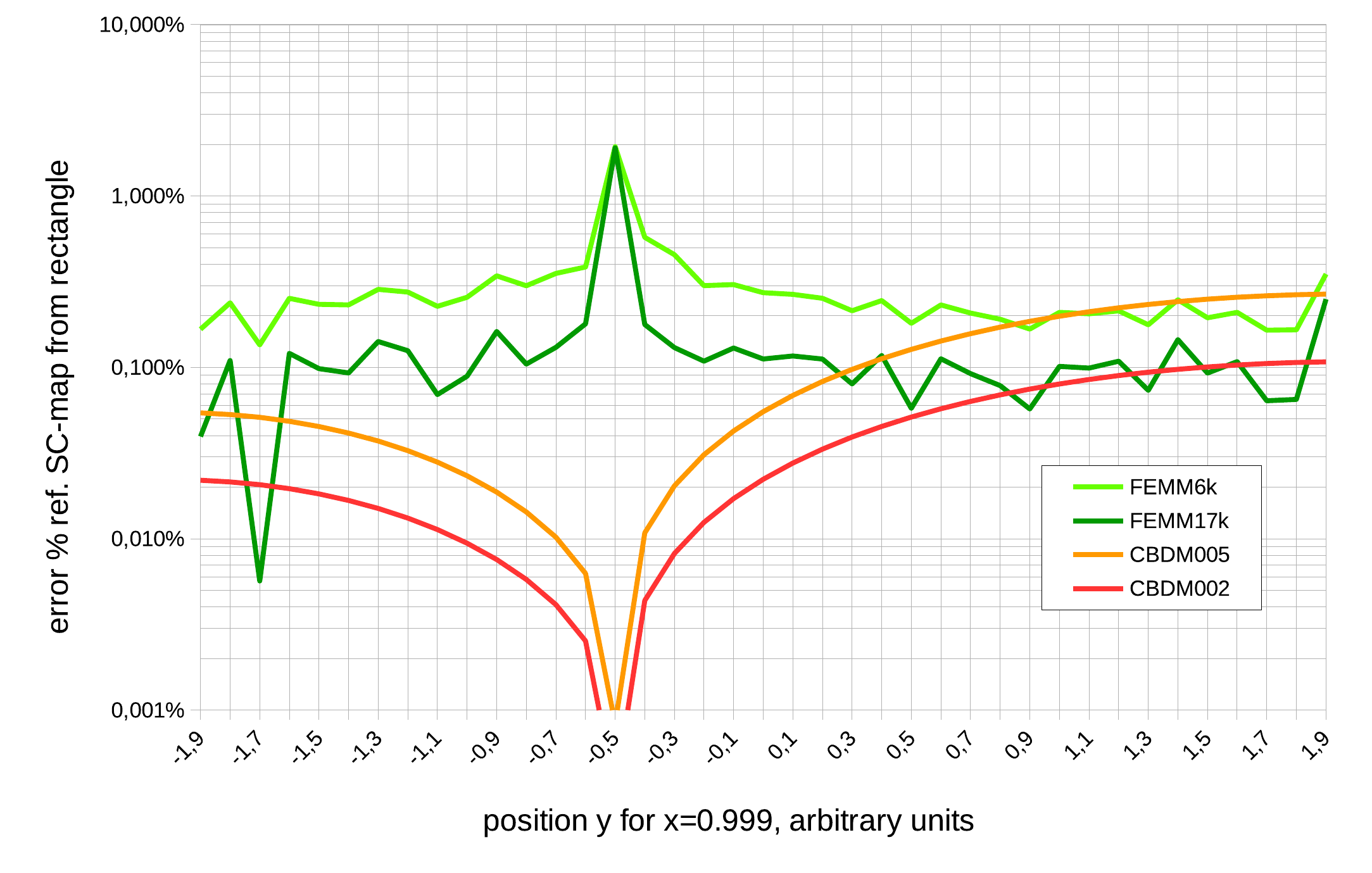}
	\caption{The error in potential distribution along the vertical line with coordinate $x=0.999$ units, without splitting $\mathbb{P}$ for the CBDM.}
	\label{fig:hom1-V}
\end{figure*}
As appears from figure \ref{fig:hom1-V}, even with a rather coarse discretization of $\Gamma_N$ the CBDM behaves as well as FEMM with high number of nodes all along the line. In particular, the results are exceptionally good right nearby the end of the slit, where the former retains the distinctive property of conformal mappings of being insensitive by nature to such kind of singularities. 

\begin{figure*}[!t]
	\centering
	\includegraphics[width=0.8\linewidth]{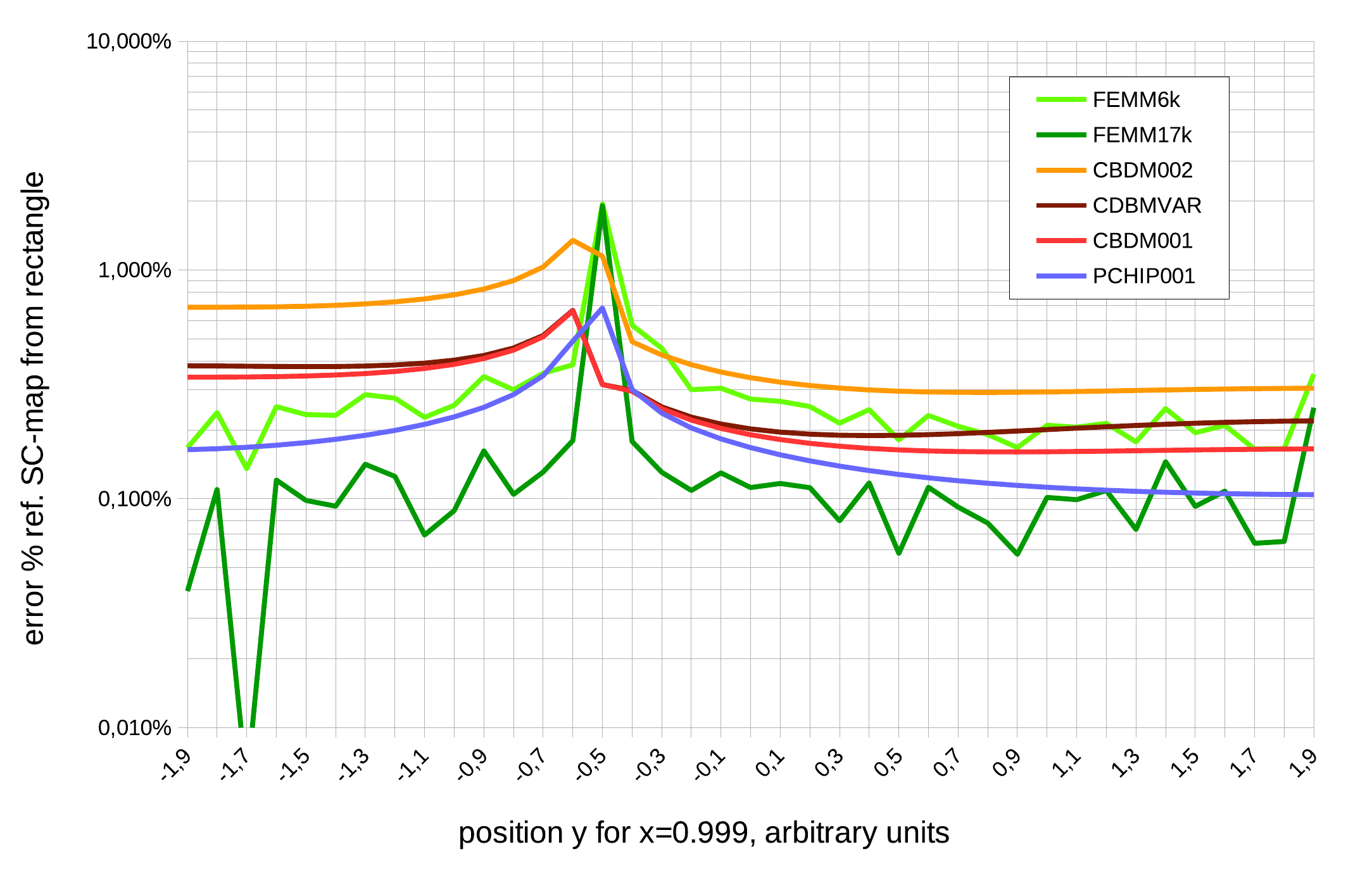}
	\caption{The error in potential distribution along the vertical line with coordinate $x=0.999$ units, after splitting $\mathbb{P}$ for the CBDM.}
	\label{fig:hom2-V}
\end{figure*}

Second, we split $\mathbb{P}$ into two subdomains as in figure \ref{fig:scgeom}-right, and compare potentials in the same way: refer now to table \ref{tab:2} and figure \ref{fig:hom2-V}.
\begin{table}[h]
	\centering
	\begin{tabular}{|c|c|p{4cm}|} \hline
		\textbf{Label}	& \textbf{Method}  & \textbf{Details}  \\ \hline\hline
		FEMM3k	& FEMM & coarser mesh, see fig. \ref{fig:FEMM_GRID}-left \\ 
		FEMM17k	& FEMM & finer mesh, see fig. \ref{fig:FEMM_GRID}-right \\ 
		CBDM002	& CBDM & step of 0.02 units on $\Gamma_N$; 293 SOR iterations \\
		CBDMVAR	& CBDM & different steps on $\Gamma_N$: 0.02 on Neumann side of $\mathbb{A}$, 0.05 on Neumann side of $\mathbb{B}$, 0.01 on $\mathbb{AB}$ interface; 393 SOR iterations \\
		CBDM001	& CBDM & step of 0.01 units on $\Gamma_N$, 595 SOR iterations \\
		PCHIP001	& CBDM & same as CDBM001, with \texttt{pchip()} interpolation of potential; 596 SOR iterations \\ \hline
	\end{tabular}
	\caption{Description of Methods in Figure \ref{fig:hom2-V}}
	\label{tab:2}
\end{table}
In short, domain partitioning breaks the magic of conformal mapping! This should not be surprising: its strengths lies on its capability of handling a domain \emph{as a whole}, and we're breaching this very point; nonetheless, results tell us that it isn't really lagging behind high density mesh FEA, and maintains the upper hand nearby the critical point indeed. One must also bear in mind that the number of points/stencils for the CBDM when using e.g. an uniform discretization step of 0.01 is only 1267/555: of course there is no direct relationship with the number of nodes/elements of figure \ref{fig:FEMM_GRID}, but admittedly this sounds as a dramatic reduction in computational effort. The aforementioned leads us to investigate whether one could obtain better results by better integration of (\ref{point:stepA}) while keeping the same step of 0.01: so we make use of the \texttt{pchip()} routine and build the Piecewise Cubic Hermite Interpolating Polynomials from $\{\psi(\theta'_k)\}_m$, integrate symbolically and evaluate between 0 and $2\pi$. The outcome, labeled PCHIP001, is an even closer approach to FEMM17k with again better behavior at the strip end; on the other hand, whereas CBDM001 with (\ref{eqn:potcbdm}) takes about 1.1 s for its solution\footnote{Time on \textsc{Matlab}; running on Octave with no JIT can take 10-15 times longer.}, the exploitation of the canned, general purpose \texttt{pchip()} as is runs 120 times slower.

As for capacitance  (table \ref{tab:3}), the map from the rectangle gives $C=1.9717349\cdot10^{-11}$ F/m for the half device, and this value is assumed exact. Both with FEMM and CBDM it is calculated as the total flux $\varepsilon\int(\nabla\psi\cdot\hat{n})$ entering the shield; for the CBDM in particular the normal component of potential gradient, entering the surface at a distance of 0.025 units from the shield\footnote{This value is so chosen as to avoid numerical difficulties with the calculation of image points in the $t$-plane.}, is calculated via (\ref{eqn:gradcbdm}) and (\ref{eqn:grad}):
\begin{table}[h]
	\centering
	\begin{tabular}{|c|c|c|} \hline
		\textbf{Label} & \textbf{Partitions of $\mathbb{P}$}  & \textbf{Error $\cdot10^{-2}$}  \\
		\hline\hline
		FEMM3k	&	& +0.28 \\
		FEMM17k	&	& +0.10 \\
		CBDM005	& none	& +0.04 \\
		CBDM002	& none	& +0.03 \\
		CBDM002	& 1	& +0.36 \\
		CBDMVAR	& 1	& +0.07 \\
		CBDM001	& 1	& +0.09 \\
		PCHIP001	& 1	& -0.06 \\ \hline
	\end{tabular}
	\caption{Error on Capacitance, ref. \texttt{rectmap}}
	\label{tab:3}
\end{table} 
the latter with no domain partitions is the big winner as expected; CBDM002 with domain partition is not too far from FEMM3k, and CBDMVAR, CBDM001 and PCHIP001 stand on par with FEMM17k.
\\

\begin{figure}[!t]
	\centering
	\includegraphics[width=1\linewidth]{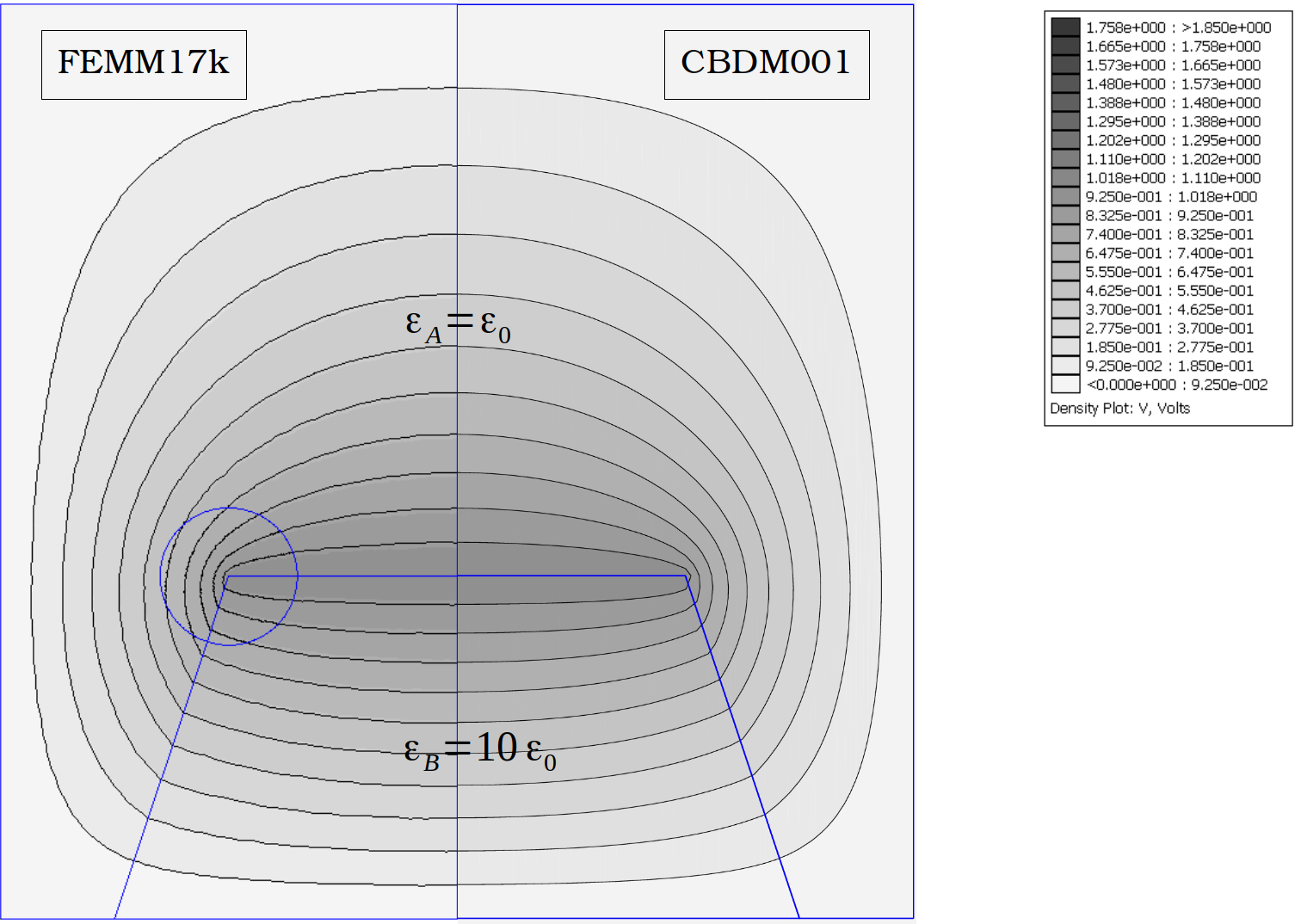}
	\caption{The whole device cross-section with the thin microstrip laying on a dielectric support of permittivity greater than air; equipotential lines are shown, as plotted and shaded by FEMM17k (left half) and CBDM001 (right half). Legend by FEMM.}
	\label{fig:field_all}
\end{figure}

When it comes to inhomogeneous domain, in our example $\varepsilon_B/\varepsilon_A=10$, we lose the support and precision of straight conformal mapping from the rectangle; so we consider the result obtained from FEMM17k the best guess, and compare with it. Figure \ref{fig:field_all} is a qualitative yet significant demonstration of the capabilities of the CBDM, obtained by juxtaposing the equipotential contour plot from FEMM17k and the horizontal-mirrored one from CBDM001.
In this case, FEMM17k returns a capacitance, calculated as before, of $C=8.46599\cdot10^{-11}$ F/m, whereas CBDM001 gives $C=8.34019\cdot10^{-11}$ F/m and PCHIP001 $C=8.32917\cdot10^{-11}$ F/m, the relative differences being of $-1.49\%$ and $-1.62\%$ respectively; these are apparently much larger than in the case of homogeneous domain, and can be a point for future investigation.

\section{Discussion}
The Conformal Boundary Differences Method (CBDM) has been introduced for the numerical solution of potential field problems. Whilst built on well-established concepts and methods of potential theory, conformal mapping and discretization by means of finite differences, it features novel aspects in how all these are made work together, also in presence of multiply connected and inhomogeneous domains, getting the Schwarz-Christoffel transformation out of its usual ancillary role with respect to FDM, FEM, BEM and the like.

This substantial degree of innovation was the reason for us to leave out the study of field sources and consider Laplacian electrostatic fields only: thus we have been able to focus on the elementary theoretical foundations, and leave larger room for the detailed discussion of a classical model problem and its comparison with results from FEA. We have found that the CBDM with a barely adequate boundary discretization is capable to stay on par with the FEM on finely meshed domains, and, being built on the SC-Toolbox, a few hundred lines of loosely optimized \textsc{Matlab}/Octave code provide accurate solutions via SOR iterations in seconds or less on a reasonably modern PC. It retains much of the strengths of conformal mappings, while broadening its scope to problems never tackled before with that technique alone. And it seems to lend itself well to the development of a consistent programming framework.

Clearly more and more case studies are needed in order to assess the merits and limits of the CBDM: here we have just scratched the surface. Laplacian fields of different nature (e.g. magnetic or thermal) are dealt with in a straightforward analogous manner; future works shall consider the presence of field sources, other types of boundary conditions (Robin, periodic) and unbounded regions. An important topic is also the analysis of errors arising from discretization and from computation, and how they affect results. This would require a dedicated study of its own and remains out of the scope of the present work, but we can rest assured of one point: as far as the potential distribution at the boundaries are adequately described and processed, the combination of the Schwarz formula and the Schwarz-Christoffel transformation can provide a solution exact up to machine precision at any point of the problem domain.

\fancyhf{}
\fancyhead[LE,RO]{\thepage}
\fancyhead[RE]{\textsc{Stefano Costa\hfill The CBDM -}\hspace{14pt}}
\fancyhead[LO]{\textsc{Stefano Costa\hfill The CBDM -}\hspace{14pt}}

\section*{References}
\bibliography{cbdm_biblio}
\bibliographystyle{alpha}
%\nocite{G51}

\begin{onecolumn}
\section*{Appendix: Working $\nabla\psi(0)$ out with Maxima}
%% Created with wxMaxima 15.08.2
\noindent
%%%%%%%%%%%%%%%
%%% INPUT:
\begin{minipage}[t]{8ex}\color{red}\bf
	\begin{verbatim}
	-->  
	\end{verbatim}
\end{minipage}
\begin{minipage}[t]{\textwidth}\color{blue}
	\begin{verbatim}
	f(x,y):=-%i*log((%e^(%i*a)-x-%i*y)/(%e^(%i*b)-x-%i*y));
	\end{verbatim}
\end{minipage}
%%% OUTPUT:

\[\displaystyle
\parbox{10ex}{$\color{labelcolor}\mathrm{\tt (\%o1) }\quad $}
\mathrm{f}\left( x,y\right) :=-i\cdot \mathrm{log}\left( \frac{-i\cdot y-x+{{e}^{i\cdot a}}}{-i\cdot y-x+{{e}^{i\cdot b}}}\right) \mbox{}
\]
%%%%%%%%%%%%%%%

\noindent
%%%%%%%%%%%%%%%
%%% INPUT:
\begin{minipage}[t]{8ex}\color{red}\bf
	\begin{verbatim}
	-->  
	\end{verbatim}
\end{minipage}
\begin{minipage}[t]{\textwidth}\color{blue}
	\begin{verbatim}
	g(x,y):=realpart(f(x,y));
	\end{verbatim}
\end{minipage}
%%%%%%%%%%%%%%%

\noindent
%%%%%%%%%%%%%%%
%%% INPUT:
\begin{minipage}[t]{8ex}\color{red}\bf
	\begin{verbatim}
	-->  
	\end{verbatim}
\end{minipage}
\begin{minipage}[t]{\textwidth}\color{blue}
	\begin{verbatim}
	load(vect);
	\end{verbatim}
\end{minipage}
%%% OUTPUT:

\[\displaystyle
\parbox{10ex}{$\color{labelcolor}\mathrm{\tt (\%o5) }\quad $}
\mbox{}
\]
/usr/share/maxima/5.37.2/share/vector/vect.mac

\noindent
%%%%%%%%%%%%%%%
%%% INPUT:
\begin{minipage}[t]{8ex}\color{red}\bf
	\begin{verbatim}
	-->  
	\end{verbatim}
\end{minipage}
\begin{minipage}[t]{\textwidth}\color{blue}
	\begin{verbatim}
	scalefactors([x,y]);
	\end{verbatim}
\end{minipage}
%%% OUTPUT:

%\[\displaystyle
%\parbox{10ex}{$\color{labelcolor}\mathrm{\tt (\%o6) }\quad $}
%\mathit{done}\mbox{}
%\]
%%%%%%%%%%%%%%%

\noindent
%%%%%%%%%%%%%%%
%%% INPUT:
\begin{minipage}[t]{8ex}\color{red}\bf
	\begin{verbatim}
	-->  
	\end{verbatim}
\end{minipage}
\begin{minipage}[t]{\textwidth}\color{blue}
	\begin{verbatim}
	gdg: grad(g(x,y));
	\end{verbatim}
\end{minipage}
%%% OUTPUT:

\begin{align*}
\displaystyle
\parbox{10ex}{$\color{labelcolor}\mathrm{\tt (\%o7) }\quad $}
\mathrm{grad}(\mathrm{atan2}\left( \frac{y-\mathrm{sin}\left( b\right) }{\sqrt{{{\left( y-\mathrm{sin}\left( b\right) \right) }^{2}}+{{\left( \mathrm{cos}\left( b\right) -x\right) }^{2}}}},\frac{\mathrm{cos}\left( b\right) -x}{\sqrt{{{\left( y-\mathrm{sin}\left( b\right) \right) }^{2}}+{{\left( \mathrm{cos}\left( b\right) -x\right) }^{2}}}}\right) - \\
-\mathrm{atan2}\left( y-\mathrm{sin}\left( a\right) ,\mathrm{cos}\left( a\right) -x\right) )\mbox{}
\end{align*}
%%%%%%%%%%%%%%%

\noindent
%%%%%%%%%%%%%%%
%%% INPUT:
\begin{minipage}[t]{8ex}\color{red}\bf
	\begin{verbatim}
	-->  
	\end{verbatim}
\end{minipage}
\begin{minipage}[t]{\textwidth}\color{blue}
	\begin{verbatim}
	ev(express(gdg), diff);
	\end{verbatim}
\end{minipage}
%%% OUTPUT:

\begin{align*}
\displaystyle
\parbox{10ex}{$\color{labelcolor}\mathrm{\tt (\%o8) }\quad $}
[\frac{{{\left( \mathrm{cos}\left( b\right) -x\right) }^{2}}\cdot \left( y-\mathrm{sin}\left( b\right) \right) }{{{\left( {{\left( \mathrm{cos}\left( b\right) -x\right) }^{2}}+{{\left( y-\mathrm{sin}\left( b\right) \right) }^{2}}\right) }^{2}}\cdot \left( \frac{{{\left( y-\mathrm{sin}\left( b\right) \right) }^{2}}}{{{\left( y-\mathrm{sin}\left( b\right) \right) }^{2}}+{{\left( \mathrm{cos}\left( b\right) -x\right) }^{2}}}+\frac{{{\left( \mathrm{cos}\left( b\right) -x\right) }^{2}}}{{{\left( y-\mathrm{sin}\left( b\right) \right) }^{2}}+{{\left( \mathrm{cos}\left( b\right) -x\right) }^{2}}}\right) }- \\
-\frac{\left( y-\mathrm{sin}\left( b\right) \right) \cdot \left( \frac{{{\left( \mathrm{cos}\left( b\right) -x\right) }^{2}}}{{{\left( {{\left( \mathrm{cos}\left( b\right) -x\right) }^{2}}+{{\left( y-\mathrm{sin}\left( b\right) \right) }^{2}}\right) }^{\frac{3}{2}}}}-\frac{1}{\sqrt{{{\left( y-\mathrm{sin}\left( b\right) \right) }^{2}}+{{\left( \mathrm{cos}\left( b\right) -x\right) }^{2}}}}\right) }{\sqrt{{{\left( y-\mathrm{sin}\left( b\right) \right) }^{2}}+{{\left( \mathrm{cos}\left( b\right) -x\right) }^{2}}}\cdot \left( \frac{{{\left( y-\mathrm{sin}\left( b\right) \right) }^{2}}}{{{\left( y-\mathrm{sin}\left( b\right) \right) }^{2}}+{{\left( \mathrm{cos}\left( b\right) -x\right) }^{2}}}+\frac{{{\left( \mathrm{cos}\left( b\right) -x\right) }^{2}}}{{{\left( y-\mathrm{sin}\left( b\right) \right) }^{2}}+{{\left( \mathrm{cos}\left( b\right) -x\right) }^{2}}}\right) }- \\
-\frac{y-\mathrm{sin}\left( a\right) }{{{\left( y-\mathrm{sin}\left( a\right) \right) }^{2}}+{{\left( \mathrm{cos}\left( a\right) -x\right) }^{2}}}, \\
\frac{\left( \mathrm{cos}\left( b\right) -x\right) \cdot {{\left( y-\mathrm{sin}\left( b\right) \right) }^{2}}}{{{\left( {{\left( \mathrm{cos}\left( b\right) -x\right) }^{2}}+{{\left( y-\mathrm{sin}\left( b\right) \right) }^{2}}\right) }^{2}}\cdot \left( \frac{{{\left( y-\mathrm{sin}\left( b\right) \right) }^{2}}}{{{\left( y-\mathrm{sin}\left( b\right) \right) }^{2}}+{{\left( \mathrm{cos}\left( b\right) -x\right) }^{2}}}+\frac{{{\left( \mathrm{cos}\left( b\right) -x\right) }^{2}}}{{{\left( y-\mathrm{sin}\left( b\right) \right) }^{2}}+{{\left( \mathrm{cos}\left( b\right) -x\right) }^{2}}}\right) }+ \\
+\frac{\left( \mathrm{cos}\left( b\right) -x\right) \cdot \left( \frac{1}{\sqrt{{{\left( y-\mathrm{sin}\left( b\right) \right) }^{2}}+{{\left( \mathrm{cos}\left( b\right) -x\right) }^{2}}}}-\frac{{{\left( y-\mathrm{sin}\left( b\right) \right) }^{2}}}{{{\left( {{\left( \mathrm{cos}\left( b\right) -x\right) }^{2}}+{{\left( y-\mathrm{sin}\left( b\right) \right) }^{2}}\right) }^{\frac{3}{2}}}}\right) }{\sqrt{{{\left( y-\mathrm{sin}\left( b\right) \right) }^{2}}+{{\left( \mathrm{cos}\left( b\right) -x\right) }^{2}}}\cdot \left( \frac{{{\left( y-\mathrm{sin}\left( b\right) \right) }^{2}}}{{{\left( y-\mathrm{sin}\left( b\right) \right) }^{2}}+{{\left( \mathrm{cos}\left( b\right) -x\right) }^{2}}}+\frac{{{\left( \mathrm{cos}\left( b\right) -x\right) }^{2}}}{{{\left( y-\mathrm{sin}\left( b\right) \right) }^{2}}+{{\left( \mathrm{cos}\left( b\right) -x\right) }^{2}}}\right) }- \\
-\frac{\mathrm{cos}\left( a\right) -x}{{{\left( y-\mathrm{sin}\left( a\right) \right) }^{2}}+{{\left( \mathrm{cos}\left( a\right) -x\right) }^{2}}}]\mbox{}
\end{align*}
%%%%%%%%%%%%%%%

\noindent
%%%%%%%%%%%%%%%
%%% INPUT:
\begin{minipage}[t]{8ex}\color{red}\bf
	\begin{verbatim}
	-->  
	\end{verbatim}
\end{minipage}
\begin{minipage}[t]{\textwidth}\color{blue}
	\begin{verbatim}
	define(gdg(x,y), %);
	\end{verbatim}
\end{minipage}
%%%%%%%%%%%%%%%

\noindent
%%%%%%%%%%%%%%%
%%% INPUT:
\begin{minipage}[t]{8ex}\color{red}\bf
	\begin{verbatim}
	-->  
	\end{verbatim}
\end{minipage}
\begin{minipage}[t]{\textwidth}\color{blue}
	\begin{verbatim}
	gdg(0,0);
	\end{verbatim}
\end{minipage}
%%%%%%%%%%%%%%%

\noindent
%%%%%%%%%%%%%%%
%%% INPUT:
\begin{minipage}[t]{8ex}\color{red}\bf
	\begin{verbatim}
	-->  
	\end{verbatim}
\end{minipage}
\begin{minipage}[t]{\textwidth}\color{blue}
	\begin{verbatim}
	trigsimp(%);
	\end{verbatim}
\end{minipage}
%%% OUTPUT:

\[\displaystyle
\parbox{10ex}{$\color{labelcolor}\mathrm{\tt (\%o11) }\quad $}
[\mathrm{sin}\left( a\right) -\mathrm{sin}\left( b\right) ,\mathrm{cos}\left( b\right) -\mathrm{cos}\left( a\right) ]\mbox{}
\]
%%%%%%%%%%%%%%%
%% End WxMaxima
\end{onecolumn}

\end{document}